\documentclass[onecolumn]{article}

\usepackage{amssymb}
\usepackage{amsmath}
\usepackage{graphicx}
\usepackage{epsfig}
\usepackage{float}
\usepackage{subfigure}
\usepackage{psfrag}
\usepackage[T1]{fontenc}
\usepackage{color}
\usepackage{array}
\usepackage{dsfont}
\usepackage{algorithm,algorithmic}
\usepackage{tikz}
\usepackage{listings}
\usepackage{epsfig}
\usetikzlibrary{arrows,decorations.pathmorphing,backgrounds,fit,positioning,shapes.symbols,chains,mindmap,trees,snakes}

\definecolor{bleuONERA}{RGB}{16,97,169}
\definecolor{grisONERA}{RGB}{64,64,66}
\providecommand{\red}[1]{\textcolor[rgb]{0.98,0.00,0.00}{#1}}

\providecommand{\green}[1]{\textcolor[rgb]{0.00,0.58,0.00}{#1}}

\providecommand{\code}[1]{\textbf{\texttt{#1}}} %

\providecommand{\matlab}[0]{\textsc{Matlab}~}

\providecommand{\mimo}[0]{\textbf{MIMO}~}

\providecommand{\siso}[0]{\textbf{SISO}~}

\providecommand{\ode}[0]{\textbf{ODE}~}
\providecommand{\pde}[0]{\textbf{PDE}~}
\providecommand{\lti}[0]{\textbf{LTI}~}

\providecommand{\tfirka}[0]{\textbf{TF-IRKA}~}

\providecommand{\mfsa}[0]{\textbf{$\Linf$-MFSA}~}

\providecommand{\ie}[0]{\emph{i.e.}~}
\providecommand{\etc}[0]{\emph{etc.}~}
\providecommand{\eg}[0]{\emph{e.g.}~}

\newenvironment{eq}{\everymath {\displaystyle \everymath{ }} \equation}{ \endequation} %
\providecommand{\norm}[1]{|| #1 ||} %
\providecommand{\trace}[0]{\mathbf{tr} } %
\providecommand{\nynu}[0]{n_{y}\times{n_{u}}}

\providecommand{\pare}[1]{\left(#1\right) } %
\providecommand{\x}[0]{\mathbf{x}} %
\renewcommand{\u}{\mathbf{u}} %
\providecommand{\y}[0]{\mathbf{y}} %
\providecommand{\yr}[0]{\mathbf{\hat{y}}} %
\providecommand{\Htranr}[0]{\mathbf{\hat{H}}} %
\providecommand{\Er}[0]{{\hat{E}}} %
\providecommand{\Ar}[0]{{\hat{A}}} %
\providecommand{\Br}[0]{{\hat{B}}} %
\providecommand{\Cr}[0]{{\hat{C}}} %
\providecommand{\Dr}[0]{{\hat{D}}} %

\providecommand{\Hsystem}[0]{\mathbf{\Sigma}} %
\providecommand{\Hreal}[0]{\mathcal{S}} %
\providecommand{\Htran}[0]{\mathbf{H}} %
\providecommand{\EE}[0]{{E}} %
\providecommand{\A}[0]{{A}} %
\providecommand{\B}[0]{{B}} %
\providecommand{\C}[0]{{C}} %
\providecommand{\DD}[0]{{D}} %

\providecommand{\Htwo}[0]{{\mathcal{H}_{2}}} %
\providecommand{\Hinf}[0]{{\mathcal{H}_{\infty}}} %
\providecommand{\Ltwo}[0]{{\mathcal{L}_{2}}} %
\providecommand{\Linf}[0]{{\mathcal{L}_{\infty}}} %

\providecommand{\Cplx}[0]{\mathbb{C}} %
\providecommand{\Real}[0]{\mathbb{R}} %
\makeindex             

\begin{document}

\title{Interpolation-based irrational model control design and stability analysis}
\author{Charles Poussot-Vassal, Pauline Kergus and Pierre Vuillemin
\thanks{Charles Poussot-Vassal and Pierre Vuillemin are with ONERA / DTIS, Universit\'e de Toulouse, F-31055 Toulouse, France and Pauline Kergus is with Lund University, Lund, Sweden. Contact:  \texttt{charles.poussot-vassal@onera.fr}}
}
%
%
\maketitle


\abstract{The versatility of data-driven approximation by interpolatory methods, originally settled for model approximation purpose, is illustrated in the context of linear controller design and stability analysis of irrational models. To this aim, following an academic driving example described by a linear partial differential equation, it is shown how the Loewner-based interpolation may be an essential ingredient for control design and stability analysis. More specifically, the interpolatory framework is first used to approximate the irrational model by a rational one that can be used for model-based control, and secondly, it is used for direct data-driven control design, showing equivalent results. Finally, this interpolation framework is employed for estimating the stability of the interconnection of the irrational model with a rational controller.}

\section{Introduction and problem statement}

\subsection{Motivations for interpolation as a pivotal tool and driving example}

Modelling, simulation, control and analysis of irrational models such as those described by linear Partial Differential Equations (\textbf{PDE}), are challenging tasks for many practitioners. Indeed, standard numerical tools developed for the rational function case are not tailored in the irrational setting and require dedicated attention (\eg eigenvalues, time-domain simulation...). In practice, engineers are often requested to discretise the irrational (infinite-dimensional) model before deploying all the numerical tools they dispose of. Beside being time consuming, this step may introduce numerical errors and lead to an iterative procedure between the control design and the model construction teams. Indeed, in an industrial context, the modelling - finite element approximation - and the control design tasks may be split in different teams and iterations to choose the ``good'' level of modelling might become an issue.

In this chapter, the control synthesis for irrational (infinite-dimensional) linear dynamical models is firstly considered. More specifically, the problem of synthesising a rational  control law  achieving some performances is  addressed through the lens of model interpolatory features. Second, the stability estimation of the closed-loop interconnection of the original linear dynamical irrational model with the synthesised linear rational controller is also done using interpolatory-like methods and approximation-oriented arguments. These arguments (tailored to linear systems only) are illustrated through a single driving example involving a transport equation controlled at a boundary. It is modelled by \eqref{eq:modelTime}, a linear \pde with constant coefficients (such equation set representing a first order linear transport equation, may be used to represent a simplified one dimensional wave equation in telecommunication, traffic jam, \etc), 
\begin{eq}
\begin{array}{rcll}
\frac{\partial \tilde\y(x,t)}{\partial x} + 2x\frac{\partial \tilde\y(x,t)}{\partial t} &=& 0 & \text{~~(transport equation)}\\ 
\tilde\y(x,0)&=&0& \text{~~(initial condition)} \\
\tilde \y(0,t)&=&\frac{1}{\sqrt{t}}*\tilde\u_f(0,t)& \text{~~(control input)}\\ 
\frac{\omega_0^2}{s^2+m\omega_0s+\omega_0^2}\u(0,s)&=&\u_f(0,s) &\text{~~(actuator model)},
\end{array}
\label{eq:modelTime}
\end{eq}
where $x\in[0~L]$, $L=3$ is the space variable and $\omega_0=3$ and $m=0.5$ are the input filter parameters. The scalar input of the model is $\tilde \u(0,t)$ (or $\u(0,s)$), the vertical force applied at the left boundary. Applying the Laplace transform, one obtains 
\begin{eq}
\frac{\partial  \y(x,s)}{\partial x} + 2x \pare{s   \y(x,s) - \tilde\y(x,0)} = 0,
\end{eq}
which solution can be given as $ \y(x,s) = a(s)e^{\int -2xs dx} = a(s)e^{-x^2s}$. Due to  boundary condition $\tilde \y(0,t)=\frac{1}{\sqrt{t}}*\tilde\u_f(t)$, we have $\y(0,s)=\frac{\sqrt{\pi}}{\sqrt{s}} \u_f(s)$, and consequently $a(s)=\frac{\sqrt{\pi}}{\sqrt{s}} \u_f(s)$. The transfer function from the input $\u(0,s)$ to the output $\y(x,s)$ reads
\begin{eq}
\y(x,s) = \frac{\sqrt{\pi}}{\sqrt{s}}e^{-x^2s} \frac{\omega_0^2}{s^2+m\omega_0s+\omega_0^2}\u(s) = \Htran(x,s) \u(0,s).
\label{eq:modelFreq}
\end{eq}

Relation \eqref{eq:modelFreq} links the (left boundary) input to the output through an irrational transfer function $\Htran(x,s)$ for any $x$ value\footnote{The exact time-domain solution of \eqref{eq:modelTime} along $x$ is given by $\tilde \y(x,t)=\tilde \u_f^{t-x^2}/\sqrt{t}$, where $\tilde \u_f$ is the output of the second order actuator transfer function, in response to $\u$.}. In addition, for the control design purpose, let us now consider that one single sensor is available, and is located at $x_m=1.9592$ along the $x$-axis\footnote{In the rest of the chapter, $x$ will be discretised with 50 points from 0 to $L=3$, and $x_m$ has been chosen to be located at $x(\lfloor{ 50 \times 2/3}\rfloor )=x(33)$, the 33-th element of $x$.}. The transfer from the same input $\u(0,s)$ to $\y_{x_m}(s)=\y(x_m,s)$ then reads:
\begin{eq}
\y_{x_m}(s)  =\Htran_{x_m}(s) \u(0,s) 
\label{eq:modelFreq_xm}
\end{eq}

\begin{figure}
\centering
\includegraphics[width=.45\columnwidth]{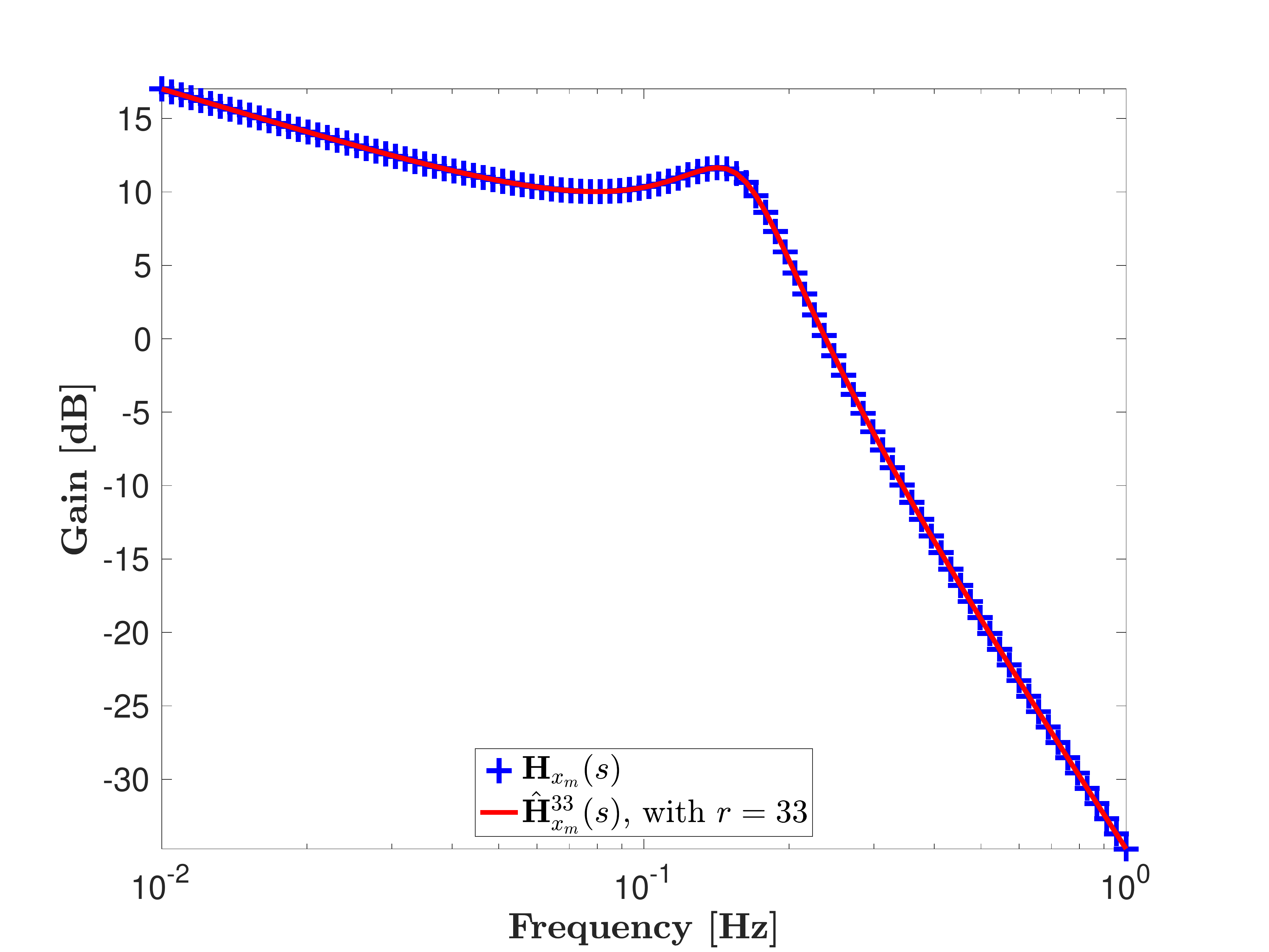}
\includegraphics[width=.45\columnwidth]{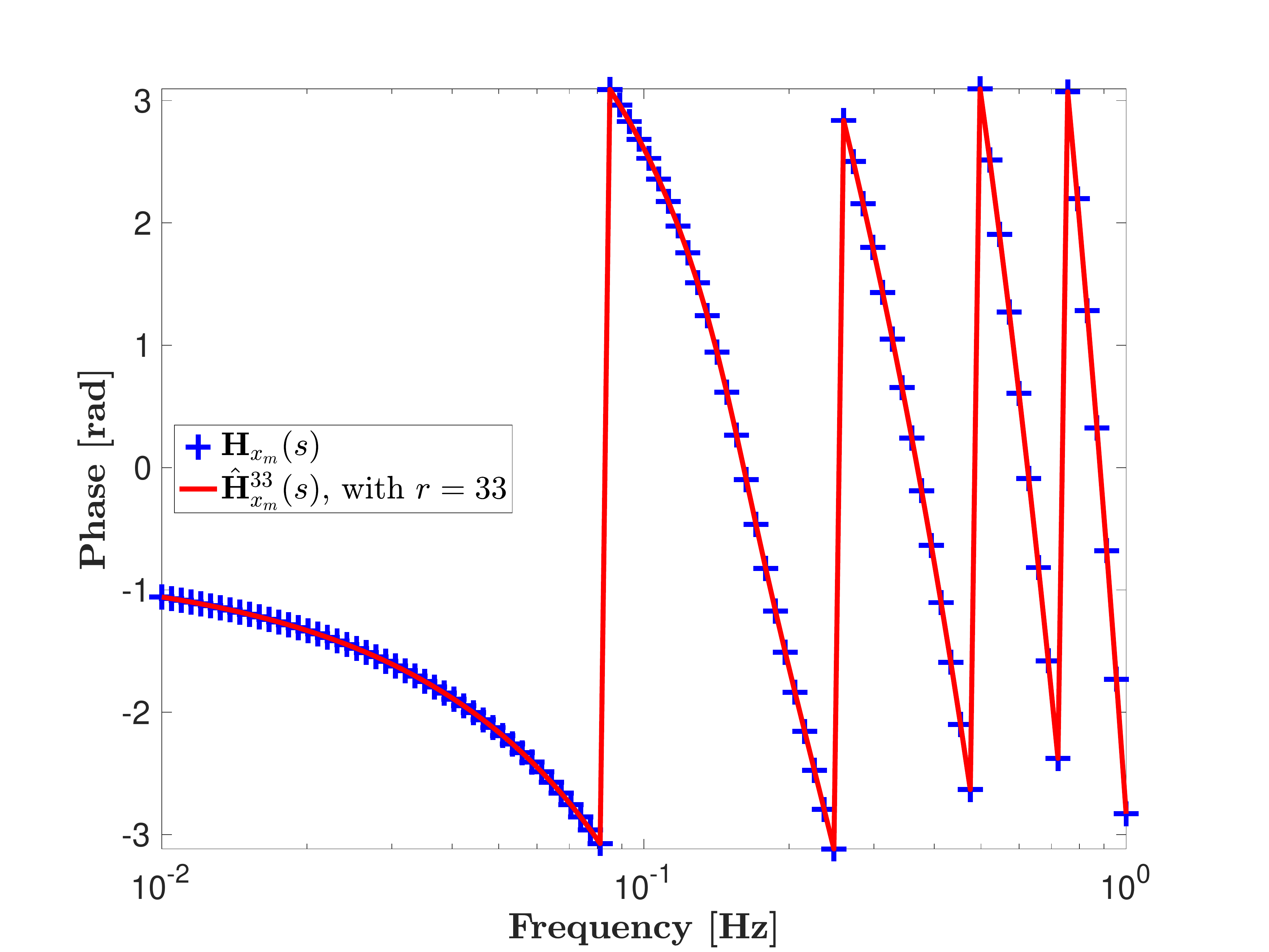}
\caption{Frequency response gain (left) and phase (right) of \eqref{eq:modelFreq_xm} (blue $+$) and its rational approximation of order $r=33$ (solid red line).}
\label{fig:modelGain}
\end{figure}

As an illustration of the above transfer function, Figure \ref{fig:modelGain} shows the frequency and phase responses at point $x=x_m$ (blue $+$). In addition, the rational approximation obtained by Loewner interpolation is also reported (solid red line).

\subsection{Problem statement, contribution and organisation}

Given a meromorphic (rational or irrational) transfer function $\Htran$, given as in \eqref{eq:modelFreq}, the objective is to design a (low order) rational feedback controller $\mathbf K$, such that the $\{\Htran,\mathbf K\}$ interconnection leads to a stable closed-loop and achieves some frequency-oriented performances. The challenge is to suggest an approach that fits to transfer function $\Htran$ being either rational or irrational. 

In this chapter, we aim at illustrating how the Loewner framework may be a pivotal tool for solving approximation problems, but also control and stability issues. 

The remainder of the chapter is organised as follows: Section \ref{sec:loe} recalls some generalities on the Loewner framework as an interpolation tool. The two proposed Loewner-driven control design methodologies are then gathered in Section \ref{sec:control}. The latter presents both the standard approximate and control approach and the data-driven control approach rooted on the Loewner framework. Both are shown to be comparable and even equivalent. Finally, in Section \ref{sec:stab}, the stability of irrational models is addressed through the lens of the Loewner tool. It is applied on the interconnection of the controllers obtained in Section \ref{sec:control} and the irrational model \eqref{eq:modelFreq}.

\subsection{Notations and preliminaries}

Along the chapter, the following notations are employed: we denote $\Htwo$ (\emph{resp.} $\Hinf$), the open subspace of $\Ltwo$ (\emph{resp.} $\Linf$) with matrix-valued function $\Htran(s)$ with $n_y$ outputs, $n_u$ inputs, $\forall s \in \Cplx$, which are analytic in $\textbf{Re}(s)> 0$ (\emph{resp.} $\textbf{Re}(s)\geq 0$). Mathematically, the $\Ltwo$ space is a vector-space of matrix valued functions defined on the imaginary axis satisfying 
$$
\int_\Real \trace\big(\overline{\Htran(\imath \omega)}\Htran(\imath \omega)^T \big)d\omega < \infty.
$$
The $\Linf$ one considers functions defined over $\Cplx_+$ satisfying 
$$
\sup_\omega \norm{\Htran(\imath \omega)}_2< \infty.
$$
Moreover $\Htwo$ and $\Hinf$ spaces consider analytic functions over the right half plane. The rational functions of the $\Hinf$ space are denoted $\mathcal{RH}_\infty$. A more detailed definition is given in the the books \cite{Partington:2004,AntoulasBook:2005}. 

Continuous \mimo \lti dynamical model (or system) $\Hsystem$ is defined as an ``input-output'' map associating an input signal $\u$ to an output one $\y$ by means of the convolution operation, defined as $\y(t) = \int_{-\infty}^{\infty} \mathbf h(t-\tau)\u(\tau)d\tau = \mathbf h(t)*\u(t)$,
where $\mathbf h(t)$ is the impulse response of the system $\Hsystem$.  It is (strictly) causal if and only if $\mathbf h(t)=0$ for ($t\leq0$) $t<0$. Then, by taking the Laplace transform of the causal convolution product above defined, one obtains 
\begin{eq}
\y(s)=\Htran(s) \u(s),
\label{eq:transfer}
\end{eq}
where $\u(s)$ and $\y(s)$ are the Laplace transform of $\u(t)$ and $\y(t)$. The $\nynu$ complex-valued matrix function $\Htran(s)$ is the transfer function of the \lti model. An \lti system $\Hsystem$ is said to be stable if and only if its transfer function $\Htran$ is bounded and analytic on $\Cplx_+$, \ie it has no singularities on the closed right half-plane. Conversely, it is said to be anti-stable if and only if its transfer function  is bounded and analytic on $\Cplx_-$ (see also \cite{Hoffman:1962} or Chapter 2 of \cite{PontesPhD:2017} for more details).

In the case where $\Htran$ is rational, it has a finite number of singularities and can be represented by a first order descriptor realisation $\Hreal:(\EE,\A,\B,\C,\DD)$ with $n_u$ inputs, $n_y$ outputs and $n$ internal variables. The model is then given by:
\begin{eq}
\EE \dot{\mathbf x}(t) =\A \x(t) +\B \u(t)  \, , \, \y(t) =\C \x(t) + \DD \u(t)
\label{eq:Hreal}
\end{eq}
where, $\x(t) \in \Real^{n}$ denotes the internal variables (the state variables if $\EE$ is invertible), and $\u(t) \in \Real^{n_u}$ and $\y(t) \in \Real^{n_y}$ are the input, output functions, respectively, while $\EE,\A \in \mathbb{R}^{n \times n}$, $\B \in \mathbb{R}^{n\times n_u}$, $\C \in \mathbb{R}^{n_y \times n}$  and  $\DD \in \mathbb{R}^{n_y \times n_u}$, are constant matrices. If the matrix pencil $(\EE,\A)$ is regular, $\Htran(s) = \C(s\EE-\A)^{-1}\B+\DD$, is called  the transfer function associated to the realisation $\Hreal$ of the system $\Hsystem$.

\section{Background on data-driven LTI model approximation}
\label{sec:loe}

Let us recall the main tool involved in this chapter, namely, the Loewner matrices. First, we define the connection between model-based and input-output data-driven approximation, then, the Loewner framework, as detailed in \cite{Mayo:2007} and \cite{AntoulasSurvey:2016} is briefly recalled. 

\subsection{\lti dynamical models and input-output data}

Given the complex-valued (rational or irrational) transfer function matrix $\Htran$ mapping the $n_u$ inputs $\u$ to the $n_y$ outputs $\y$ (as in \eqref{eq:transfer}) or the input-output data collection $\{z_i,\mathbf{\Phi}_i\}_{i=1}^N=\{z_i,\Htran(z_i)\}_{i=1}^N$  (where $z_i\in\Cplx$ and $\mathbf \Phi_i\in \Cplx^{n_y\times n_u}$) defined as,
\begin{eq}
\y(z_i) = \Htran(z_i) = \mathbf \Phi_i \u(z_i),
\label{eq:data}
\end{eq}
the approximation problem aims at constructing the (approximate) rational transfer function matrix $\Htranr$ mapping inputs $\u$ to the approximate outputs $\yr$ such that
\begin{eq}
\yr(s) = \Htranr(s) \u(s).
\label{eq:transferRed}
\end{eq}
Obviously, some objective are that \emph{(i)} the reduced inputs to outputs map should be "close" to the original, \emph{(ii)} the critical system features and structure should be preserved, and,  \emph{(iii)} the strategies for computing the reduced system should be numerically robust and stable. Approximating \eqref{eq:transfer} with \eqref{eq:transferRed} is a \emph{model-based} approximation, while, approximating \eqref{eq:data} with \eqref{eq:transferRed} belongs to the \emph{data-driven} family (see \cite{AntoulasBook:2005,PoussotHDR:2019} for examples). In this chapter, we follow the \emph{data-driven} philosophy, and more specifically the interpolation-based approach using the Loewner framework.

\subsection{Data-driven approximation and Loewner framework at a glance}

The main elements of the Loewner framework are recalled here in the single-input single-output (\textbf{SISO}) case and readers may refer to \cite{Mayo:2007} for a complete description and extension to the \mimo one. The Loewner approach is a data-driven method building a (reduced) rational descriptor \lti dynamical model $\Htran^m$ ($\Htranr^r$) of dimension $m$ ($r<m$) of the same form as \eqref{eq:Hreal}, which interpolates frequency-domain data given as \eqref{eq:data}. More specifically, let us consider a set of distinct interpolation points  $\{ z_i \}_{i=1}^{2m} \subset \Cplx$ which is split in two subsets of equal length as $\{z_i\}_{i=1}^{2m} =\{\mu_i\}_{i=1}^{m} \cup \{\lambda_i\}_{i=1}^{m}$.  The method consists in building the Loewner and shifted Loewner matrices as,
\begin{eq}
    \left [ \mathbb{L} \right ]_{ij} = \frac{\Htran(\mu_i) - \Htran(\lambda_j)}{\mu_i - \lambda_j}
    \text{ and }
        \left [ \mathbb{L}_s \right ]_{ij} = \frac{\mu_i \Htran(\mu_i) - \lambda_j \Htran(\lambda_j)}{\mu_i - \lambda_j}.
\end{eq}
The model $\Htran^m$ that interpolates $\Htran$ is given by the following descriptor realisation,
\begin{eq}
\EE^m  \dot{\mathbf{x}}(t) = \A^m \x(t) + \B^m \u(t)\, , \, \y(t) =\C^m \x(t)
\label{eq:descr}
\end{eq}
where $\EE^m = -\mathbb{L}$, $\A^m = -\mathbb{L}_s$, $[\B^m]_i = \Htran(\mu_i)$ and $[\C^m]_i = \Htran(\lambda_i)$ ($i=1,\ldots,m)$. Assuming that the number $2m$ of available data is large enough, then  it has been shown in \cite{Mayo:2007} that a minimal model $\Htranr^r$ of dimension $r < m$ that still interpolates the data\footnote{The model $\Htranr^r$ interpolates the original model $\Htran$ at the $z_i$ points. This is the reason why this method is also known as an interpolation one.} can be built with a projection of \eqref{eq:descr} provided that, for $i=1,\ldots,2m$ (note that to avoid complex arithmetic, $m$ points and their conjugate are selected),
\begin{equation}
    rank(z_i \mathbb{L} - \mathbb{L}_s) = rank( [\mathbb{L}\,\, \mathbb{L}_s]) = rank([\mathbb{L}^T\,\, \mathbb{L}_s^T]^T) = r.
    \label{eq:rankCond}
\end{equation}
In that case, denoting $Y \in \Cplx^{m \times r}$ and (\emph{resp.} $X \in \Cplx^{m \times r}$) the matrix containing the first $r$ left (\emph{resp.} right) singular vectors of $[\mathbb{L}\,\, \mathbb{L}_s]$ (\emph{resp.} $[\mathbb{L}^T\,\, \mathbb{L}_s^T]^T$). Then, $\Er^r = Y^H \EE^m X$, $\Ar^r = Y^H \A^m X$,  $\Br^r = Y^H \B^m$, and $\Cr^r = \C^m X$, is a realisation of this model $\Htranr^r$ with a McMillan degree equal to $rank(\mathbb{L})$. Note that if $r$ in \eqref{eq:rankCond} is superior to $rank(\mathbb{L})$ then $\Htranr^r$ can either have a direct feed-through $\Dr\neq 0$ or a polynomial part. In the rest of the paper, one assumes that no polynomial term is present in the state-space realisation \eqref{eq:descr}. 

\section{Model- and data-driven control synthesis paradigm}
\label{sec:control}

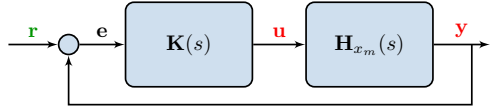
\begin{figure}[h]
\centering
\scalebox{.8}{
\tikzstyle{block} = [draw, thick,fill=bleuONERA!20, rectangle, minimum height=4em, minimum width=6em,rounded corners]
\tikzstyle{sum} = [draw, thick,fill=bleuONERA!20, circle, node distance=1cm]
\tikzstyle{input} = [coordinate]
\tikzstyle{output} = [coordinate]
\tikzstyle{pinstyle} = [pin edge={to-,thick,black}]
\tikzstyle{connector} = [->,thick]
\begin{tikzpicture}[auto, node distance=2cm,>=latex']
    \node [input, name=input] {};
    \node [sum, right of=input] (sum) {};
    \node [block, right of=sum] (controller) {$\mathbf K(s)$};
    \node [block, right of=controller, node distance=3cm] (system) {$\Htran_{x_m}(s)$};
    \draw [connector] (controller) -- node[name=u] {$\red{\u}$} (system);
    \node [output, right of=system] (output) {};
    \draw [connector] (input) -- node {$\green{\mathbf r}$} (sum);
    \draw [connector] (sum) -- node {$\mathbf e$} (controller);
    \draw [connector] (system) -- node [name=y] {$\red{\y}$}(output);
    \draw [connector] (output)+(-0.3cm,0) -- ++(-0.3cm,-1cm) -| node [near start] {} (sum.south);
\end{tikzpicture}
}
\caption{Feedback loop scheme. The objective is to find $\mathbf K(s)$ stabilising $\Htran_{x_m}(s)$ and achieving some performances.}
\label{fig:feedbackInitial}
\end{figure}

Based on the above tools, let us now jump in the first contribution of the chapter, namely, the design of a feedback controller for irrational models. This objective, presented in Figure \ref{fig:feedbackInitial}, aims at seeking for a controller $\mathbf K \in \mathcal{RH}_\infty$ such that the interconnection closed-loop $\{\Htran,\mathbf K\}$ is stable and achieves some performance \eg minimises some $\Hinf$-norm or track some closed-loop performances. 

In the rest of this section, two control design approaches are developed. The first one is a - standard - approximate and model-driven method, while the second is a direct data-driven one. The model-driven control design is based on a rational approximation of model \eqref{eq:modelFreq} using the Loewner framework, followed with a structured $\Hinf$-norm oriented control design step (note that any other approach can be used), while the data-driven relies on the Loewner framework to directly identify the controller, on the basis of the data only.

\subsection{Model-driven approximation and control}
\label{sec-model-based}

On Figure \ref{fig:modelGain}, the transfer from the boundary input $\u$ to the measurement point $\y_{x_m}$ of both the original model $\Htran_{x_m}$ and of the rational approximation one $\Htranr_{x_m}^{33}$, were illustrated. $\Htranr_{x_m}^{33}$ (of dimension $r=33$) has been obtained by sampling $\Htran_{x_m}$ as
$\{\Htran_{x_m}(z_i)\}_{i=1}^{2m}$
where $z_i=\imath\omega_i$ such that $z_i$ are closed under conjugation and for $\omega_i$ logarithmically spaced between $2\pi 10^{-2}$ and $2\pi$ ($m=200$).  At this point, thanks to the Loewner approach, it is both possible to simulate the equation using a standard \ode solver and to design a control law using rational model-based control design methods. 

\begin{figure}
\centering
\includegraphics[width=.9\columnwidth]{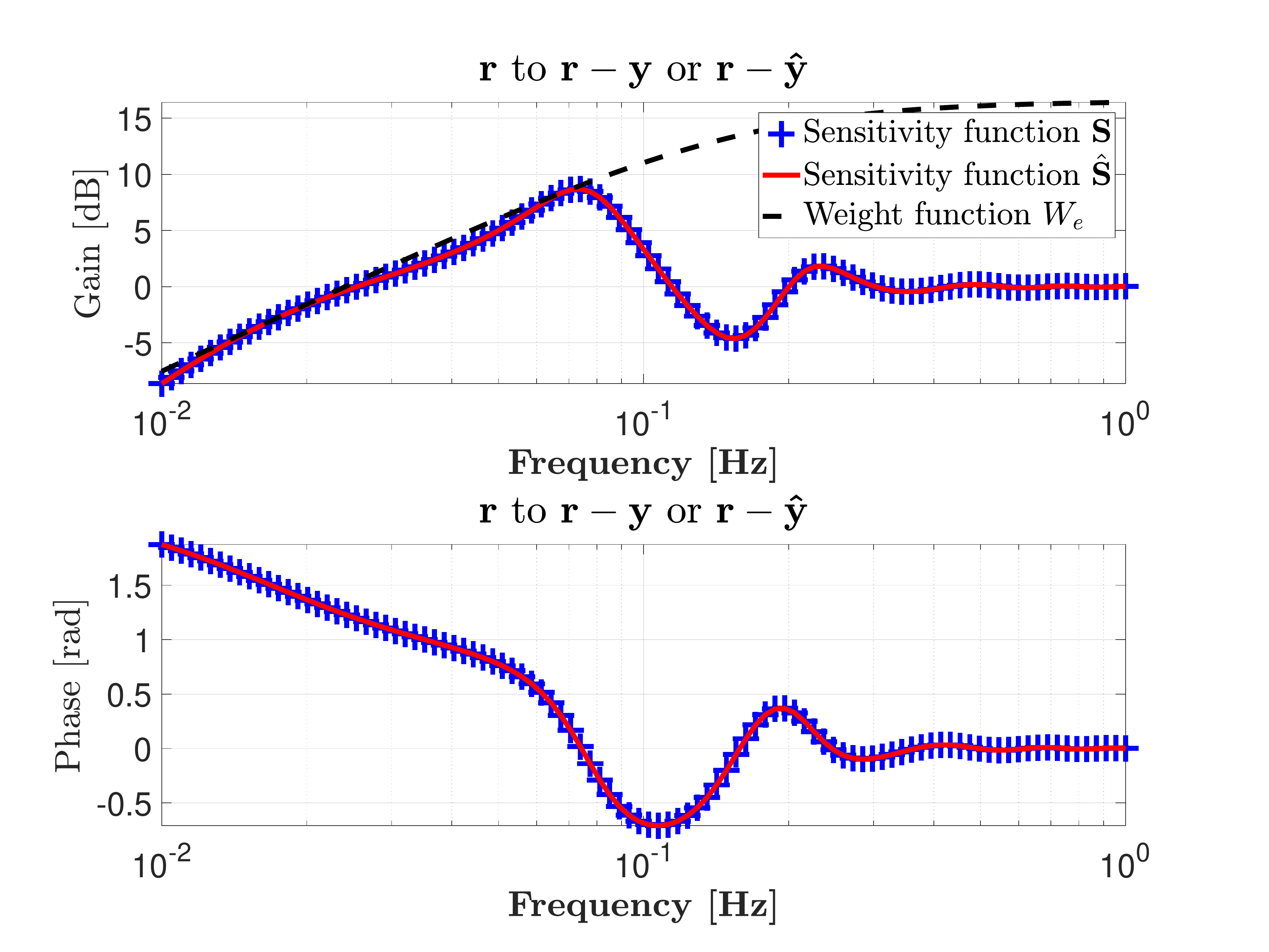}
\caption{Sensitivity function $\mathbf S$ (blue $+$)  and $\mathbf{\hat S}$ (solid red). Weighting function $W_e$ used in the control design (solid black).}
\label{fig:S}
\end{figure}

\begin{figure}
\centering
\includegraphics[width=.9\columnwidth]{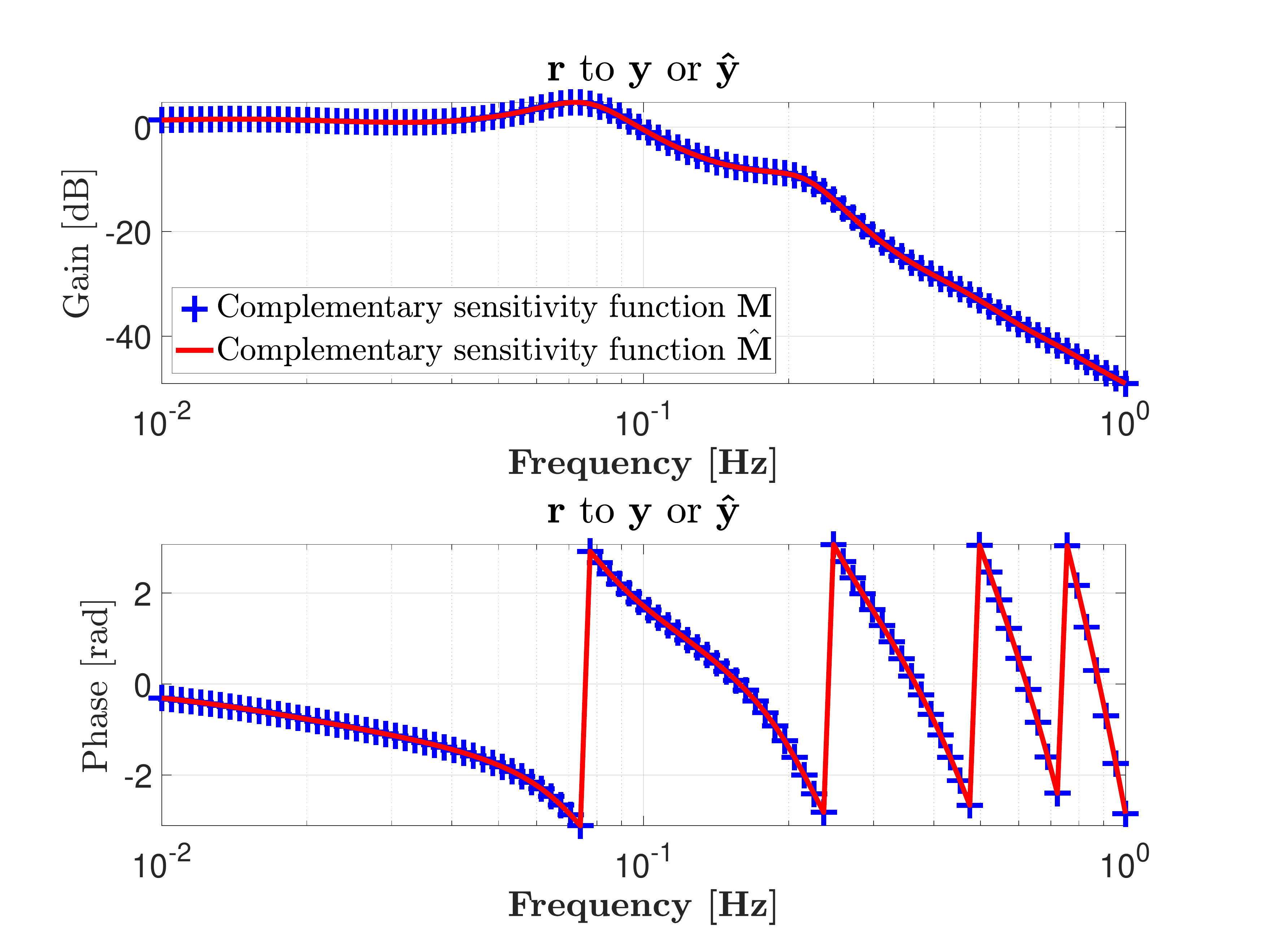}
\caption{Complementary sensitivity function $\mathbf M$ and $\mathbf{\hat M}$, linking $\mathbf r$ to $\mathbf y$ (blue $+$) or $\yr$ (solid red).}
\label{fig:M}
\end{figure}

Using the rational approximation $\Htranr_{x_m}^{33}$ of the irrational model $\Htran_{x_m}$, standard feedback synthesis methods can be applied (\ie design a feedback loop as on Figure \ref{fig:feedbackInitial}, but involving $\Htranr_{x_m}^{33}$ instead of $\Htran_{x_m}$). In this examples, the \code{hinfstruct} function embedded the \matlab Robust Control Toolbox has been used \cite{Apkarian:2006}; it allows designing fixed structure controllers while minimising some $\Hinf$-norm oriented performance transfer. Starting from $\Htranr_{x_m}^{33}$, let us first define the following generalised plant $\mathbf{T}=\Htranr_{x_m}^{33} \mathbf W_o$, where $\mathbf W_o$ is the weighting filter defining the output signals on which the $\Hinf$-norm optimisation will be done. $\mathbf W_o$ is constructed to define the desired closed-loop performances attenuation and its bandwidth. The resulting state-space realisation of the generalised plant $\mathbf T$ is then given by
\begin{eq}
\left\{
\begin{array}{rcl}
\boldsymbol{\dot{\xi}}(t) &=& \A_\xi \boldsymbol \xi(t) + \B_1 \boldsymbol r(t) + \B_2 \u(t)  \\
\mathbf z(t)&=& \C_1\boldsymbol \xi(t) + D_{11} \mathbf r(t) + D_{12} \u(t)\\
\mathbf e(t)&=& \mathbf r(t) - \y
\end{array} \right. 
\label{eq:gene}
\end{eq}
where $\boldsymbol \xi(t) \in \Real^{N}$ and $\mathbf z(t) \in \Real^{n_z}$ are the states (model plus weight state variables) and performance output signals, respectively. The performance transfer from $\mathbf r$ to $\mathbf z$, is defined as $\mathbf T_{\mathbf r\mathbf z}=\Htranr_{x_m}^{33} \mathbf W_o$. In the considered case, one aims at tracking the reference signal $\mathbf r$ and limiting the control action $\mathbf u$. One can then construct $W_o=\textbf{blkdiag}\big(W_e,W_u\big)=\textbf{blkdiag}\big(10\frac{s+1}{s},\frac{s+10}{s+1000}\big)$ describing performance  output $\mathbf z=\textbf{blkdiag}\big(W_o\mathbf e,W_u\u\big)$. The $W_o$ weighting filter has been chosen to weight the sensitivity function and guarantee no steady-state error (\eg roll-off in low frequency) and a bandwidth around $10^{-1}rad/s$. The $W_u$ one is instead used to weight actuator action in high frequencies (here the actuator will roll-off above $10$rad/s). It can be mentioned that it is also a fairly standard way of weight selection. The $\Hinf$ control design consists in finding the controller $\mathbf K$, mapping $\mathbf e$ to $\u$, such that,
\begin{eq}
\small
\mathbf K := \arg \min_{\tilde{\mathbf{K}}\in{\mathcal K}} \norm{\mathcal F_l\big(\mathbf{T}_{\mathbf r\mathbf z},\tilde{\mathbf{K}}\big)}_{\mathcal H_\infty}
\label{eq:pbHinf}
\normalsize
\end{eq}
where $\mathcal F_l(\cdot,\cdot)$ is the lower fractional operator defined as (for appropriate partitions of $M$ and $K$) by ${\cal F}_l(M,K) =  M_{11}+M_{12}K(I-M_{22}K)^{-1}M_{21}$ \cite{MagniLFR:2006}. With reference to \eqref{eq:pbHinf}, it is possible to define the class $\mathcal K$ of $\mathbf K$ to be restricted to the Proportional Integral one, meaning that one is seeking for $\mathbf K$ with the following form, $\mathbf K(s) = k_p + k_i\frac{1}{s}$, where $k_p,k_i \in \Real$. After optimisation, one obtains $k_p=0.191$ and $k_i=0.0252$ (note also that in this case, the optimal attenuation reached is $\gamma_\infty=66.954$) \footnote{The optimisation is done using the \code{hinfstruct} routine, allowing minimising the closed-loop interconnection of $\mathbf{T_{rz}}$ with $\mathbf{\tilde K}$. In general, one seek for $ \norm{\mathcal F_l\big(\mathbf{T}_{\mathbf r\mathbf z},{\mathbf{K}}\big)}_{\mathcal H_\infty}=\gamma_\infty\leq 1$. Here we simply aim at reaching stability and tracking performances.}. Figure \ref{fig:S} then shows the sensitivity function $\mathbf S$ (transfer from $\mathbf r$ to $\mathbf e$) applied both on the rational approximation $\Htranr_{x_m}^{33}$ and the original model $\Htran_{x_m}$. It shows that good tracking in low frequencies is ensured, as well as some margin properties. In addition, the complementary sensitivity function, $\mathbf M=1-\mathbf S$, is reported on Figure \ref{fig:M}, illustrates the closed-loop transfer from the reference $\mathbf r$ to the model output $\y$ or $\yr$.

The resulting control law gives similar results for both approximated (rational) and original (irrational) models thus validating the approach.

\subsection{Data-driven control}

So far, the control design has been done in a fairly standard way, involving the approximated rational model. 
Instead of designing a controller on the basis of a rational reduced-order model, as the true system's behaviour is known through \eqref{eq:transfer}, a data-driven control strategy, based on the approach presented in \cite{KergusWC:2017}, is followed. Authors stress that the main contribution of this section with respect to  \cite{KergusWC:2017}  stands in the comparison of the data-driven approach with the model-based one, resulting, as we will see later in the chapter, in exactly similar results. 

\subsubsection{Introducing the ideal controller}

Data-driven control consists in recasting the control design problem as an identification one. The main advantage of this strategy is that it provides a controller tailored to the actual system. This change of paradigm shift the identification / simplification process of model to the controller directly. Different techniques have been proposed, see references in \cite{KergusWC:2017}, considering a set of structured controllers. Recently, \cite{KergusWC:2017} pushed the interpolatory framework, and especially the Loewner one, into the process, enabling \mimo controllers design without a-priori structure selection. 

\begin{figure}[h]
\centering
\scalebox{.78}{
\tikzstyle{block} = [draw, thick,fill=bleuONERA!20, rectangle, minimum height=4em, minimum width=6em,rounded corners]
\tikzstyle{sum} = [draw, thick,fill=bleuONERA!20, circle, node distance=1cm]
\tikzstyle{input} = [coordinate]
\tikzstyle{output} = [coordinate]
\tikzstyle{pinstyle} = [pin edge={to-,thick,black}]
\tikzstyle{connector} = [->,thick]
\begin{tikzpicture}[auto, node distance=2.5cm,>=latex']
    \node [input, name=input] {};
    \node [sum, right of=input] (sum) {};
    \node [block, right of=sum] (controller) {$\mathbf K(s)$};
    \node [block, right of=controller, node distance=3cm] (system) {$\mathbf H_{x_m}(s)$};
    \node [sum, right of=system, node distance=2cm] (sum2) {};
    \node [block, above of=controller, node distance=2cm] (obj) {$\mathbf M(s)$};
    \draw [connector] (controller) -- node[name=u] {$\red{\u}$} (system);
    \node [output, right of=sum2, node distance=1cm] (output) {};
    \draw [connector] (input) -- node {$\green{\mathbf r}$} (sum);
    \draw [connector] (sum) -- node {$\mathbf e$}(controller);
    \draw [connector] (system) -- node [name=y] {$\red{\mathbf y}$}(sum2);
    \draw [connector] (sum2)+(-0.3cm,0) -- ++(-0.3cm,-1cm) -| node [near start] {} (sum.south);
    \draw [connector] (sum2) -- node {$\mathbf \varepsilon$}(output);
    \draw [connector] (input)+(0.3cm,0) |- (obj) -| (sum2);
\end{tikzpicture}
}
\caption{Data-driven control problem formulation: $\mathbf{M}$ is the reference model (objective) and $\mathbf{K}$ the controller to be designed.}
\label{fig:problem_formulation}
\end{figure}
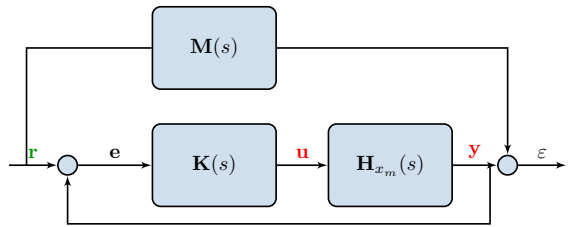

The objective is to find a controller $\mathbf{K}$ minimising the difference between the resulting closed-loop and a given reference model $\mathbf{M}$ (see Figure \ref{fig:problem_formulation}). This is made possible through the definition of the ideal controller $\mathbf{K}^\star$, being the \lti controller that would have given the desired reference model behaviour if inserted in the closed-loop. It is then defined as follows:
\begin{eq}
    \mathbf{K}^\star=\Htran_{x_m}^{-1}\mathbf{M}(I-\mathbf{M})^{-1}.
    \label{Kideal}
\end{eq}

Given a reference model $\mathbf{M}$ and provided frequency-domain data $\left\{\imath\omega_i,\mathbf{\Phi}_i\right\}_{i=1}^{2m}$ from the plant $\mathbf{H}_{x_m}$, it is possible to evaluate the frequency-response of the ideal controller $\mathbf{K}^\star$ at these very same frequencies. The main idea of the Loewner Data-Driven Control (\textbf{LDDC}) algorithm introduced in \cite{KergusWC:2017} is to interpolate the frequency-response of the ideal controller $\mathbf{K}^\star$ and to reduce it to an acceptable order.  For the present example, $m=200$ samples of the frequency-response, logarithmically spaced between $2\pi10^{-2}$ and $2\pi$ rad.s$^{-1}$ are considered (similar to the one used  for the computation  of a rational reduced-order model in Section \ref{sec:loe}). 

\subsubsection{Data-driven control design and model-free stability analysis}
As explained in \cite{kergus2019}, the reference model choice is a key factor for the \textbf{LDDC}  success, as for any other model reference control techniques. Indeed, this latter should not only represent a desirable closed-loop behaviour, but also an achievable dynamic of the considered system (\ie the ideal controller should not internally destabilise the plant and imply terrible dynamics). A reference model is said to be achievable by the plant if the corresponding ideal controller internally stabilises the plant.

The closed-loop performances are limited by the system's right hand side poles $\{p_j\}_{j=1}^{n_p}$ and zeros $\{z_i\}_{i=1}^{n_z}$ (and their respective output directions in the \mimo case) defined as $\Htran_{x_m}(z_i) = 0$ and $\Htran_{x_m}(p_j) = \infty$.
Finally, the class $\mathcal{M}$ of achievable reference models is defined as follows: 
\begin{eq}
\mathcal{M}=\big\{\mathbf{M}\in \mathcal{H}_\infty : \, \mathbf{M}(z_i) = 0 \text{ and } \mathbf{M}(p_j) = 1\big\}.
\end{eq}
In the general case, the right half plane poles and zeros of the system are estimated in \cite{kergus2019} in order to build an achievable reference model on the basis of the initial specifications given by the user. This is made possible through a data-driven stability analysis introduced in \cite{cooman2018model} and the associated estimation of instabilities presented in \cite{cooman2018estimating}. 

Here, the \pde describing the system's dynamic is known and allows to determine the performance limitations without applying the aforementioned data-driven stability analysis. The instabilities are $\Htran_{x_m}(+\infty)=0$ and $\Htran_{x_m}(0)=\infty$, implying that the reference model should satisfy $\mathbf{M}(+\infty)=0$ and $\mathbf{M}(0)=1$.

As shown in \cite{kergus2019appendix}, once the set $\mathcal{M}$ of achievable reference models is determined, the one chosen in this set influences a lot the control design. To illustrate this point, the \textbf{LDDC} algorithm is applied on the proposed example using two reference models $\mathbf{M}_1$ and $\mathbf{M}_2$. $\mathbf{M}_1$ is chosen as a perfectly damped second-order model with a natural frequency $\omega_0=0.5\textnormal{ rad.s}^{-1}$ and reads
$$
\mathbf{M}_1(s)=\frac{1}{s^2/\omega_0^2+2s/\omega_0+1},
$$
satisfying the performance limitations of the system. The second reference model $\mathbf{M}_2$ is the closed-loop obtained in the model-based control design obtained in Section \ref{sec-model-based} ($\mathbf{M}_2=\hat{\mathbf{M}}$, see Figure \ref{fig:M})\footnote{Since the procedure used to get $\hat{\mathbf{M}}$ preserves internal stability, the obtained closed-loop is necessarily achievable by the plant $\Htran_{x_m}$.}.

Once the reference models $\mathbf{M}_1$ and $\mathbf{M}_2$ are chosen, by following \eqref{Kideal}, it is possible to compute the frequency-responses of the associated ideal controllers, denoted $\mathbf{K}^\star_1$ and $\mathbf{K}^\star_2$ respectively, at the frequencies where data from the plant are available. In order to obtain a controller model $\mathbf{K}$, the Loewner framework is then applied considering the following interpolatory conditions: $\forall i=1\dots 2m$, $\mathbf{K}(\imath \omega_i)=\mathbf{K}^\star(\imath \omega_i)$.

In the present case, minimal realisations of $\mathbf{K}^\star_1$ and $\mathbf{K}^\star_2$ of order $n_1=34$ and $n_2=40$ respectively are obtained. In order to compare the results of the model-based approach with the data-driven strategy, the ideal controllers $\mathbf{K}^\star_1$ and $\mathbf{K}^\star_2$ are reduced up to a first order (using the rank revealing factorisation embedded in the Loewner framework), giving two controller models denoted $\mathbf{K}_1$ and $\mathbf{K}_2$, which transfer functions are given in \eqref{LDDC_controllers}. 
\begin{eq}
    \mathbf{K}_1(s)=\frac{0.02277}{(s+0.0382)} \text{ and } \mathbf{K}_2(s)=\frac{0.1914 (s+0.1315)}{s}=0.1914+\frac{0.0252}{s}.
    \label{LDDC_controllers}
\end{eq}

\begin{figure}[h]
\centering
\includegraphics[width=0.9\textwidth]{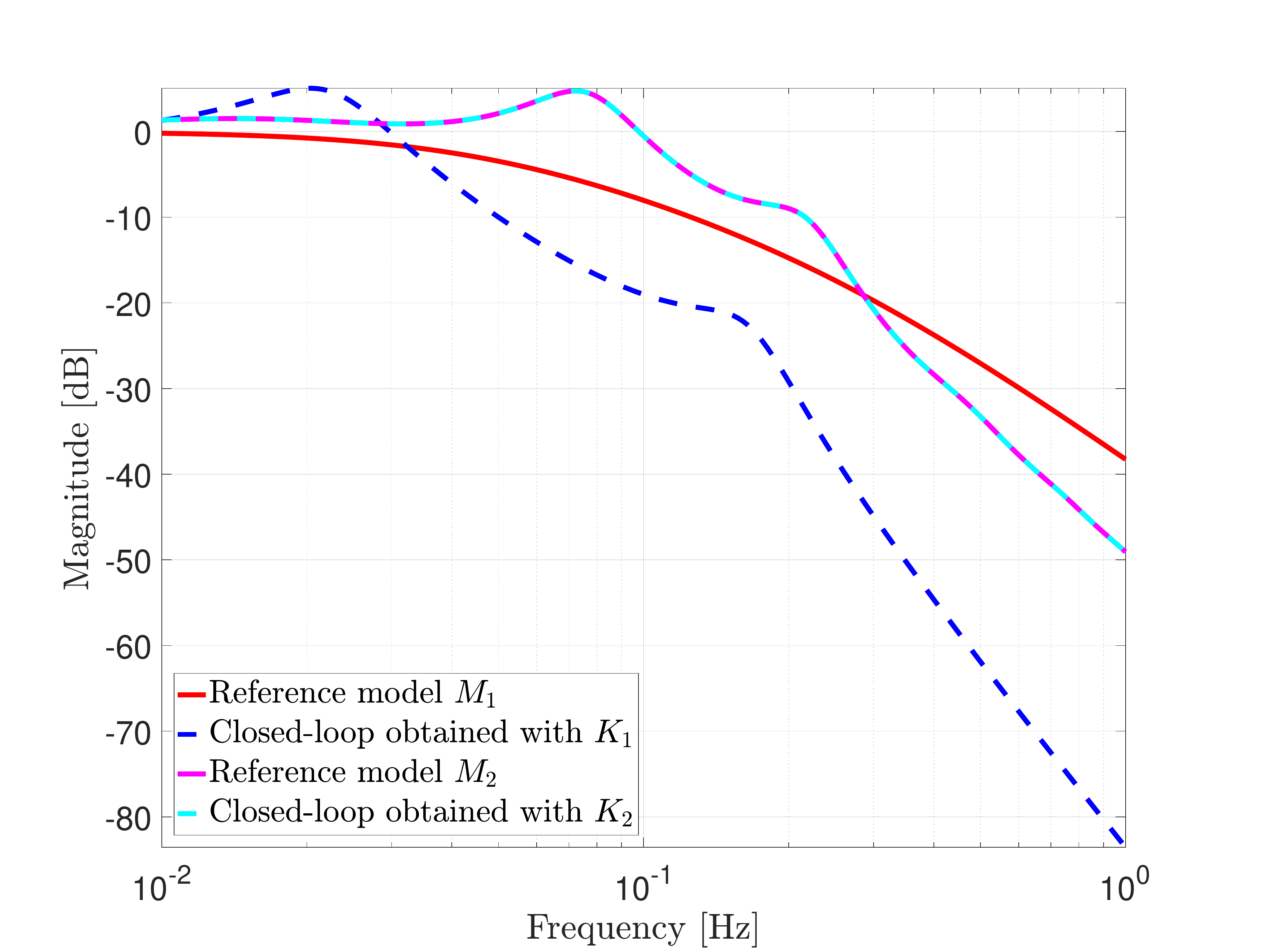}
\caption{Resulting closed-loop transfers involving $\mathbf K_1$ and $\mathbf K_2$, obtained using $\mathbf M_1$ and $\mathbf M_2$ reference models.}
\label{fig:M_DDC}
\end{figure}

When using $\mathbf{M}_2$ as reference model, the first order controller $\mathbf{K}_2$ has exactly the same expression than the one obtained in the model-based approach solving \eqref{eq:pbHinf}. The frequency responses of the resulting closed-loops obtained with $\mathbf K_1$ and $\mathbf K_2$ are visible on Figure \ref{fig:M_DDC}.

Interestingly, with reference to Figure \ref{fig:M_DDC}, $\mathbf K_2$ perfectly recovers the requested performance of $\mathbf M_2$ with a controller of rational order one (indeed, we expected to observe this result since we knew from the model-based approach presented in Section \ref{sec-model-based} that a rational control of order leading to this performance was achievable). Conversely, $\mathbf K_1$ (reduced to a rational form of dimension one) is not able to recover the performances of  $\mathbf M_1$, and a higher degree would be expected (as it is not the topic of the chapter, this point is not detailed further).

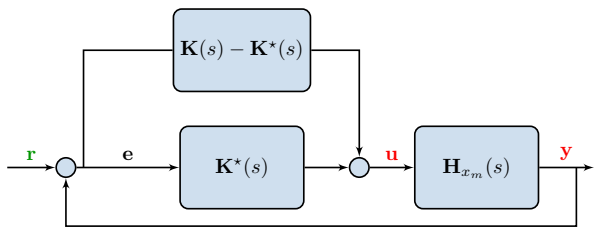
\begin{figure}[h]
\centering
\scalebox{.78}{
\tikzstyle{block} = [draw, thick,fill=bleuONERA!20, rectangle, minimum height=4em, minimum width=6em,rounded corners]
\tikzstyle{sum} = [draw, thick,fill=bleuONERA!20, circle, node distance=1cm]
\tikzstyle{input} = [coordinate]
\tikzstyle{output} = [coordinate]
\tikzstyle{pinstyle} = [pin edge={to-,thick,black}]
\tikzstyle{connector} = [->,thick]
\begin{tikzpicture}[auto, node distance=3cm,>=latex']
    \node [input, name=input] {};
    \node [sum, right of=input] (sum) {};
    \node [block, right of=sum] (controller) {$\mathbf K^\star(s)$};
    \node [sum, right of=controller, node distance=2cm] (sum2) {};
    \node [block, right of=controller, node distance=4cm] (system) {$\mathbf H_{x_m}(s)$};
    \node [block, above of=controller, node distance=2cm] (err) {$\mathbf K (s) - \mathbf K ^\star (s)$};
    \draw [connector] (controller) -- (sum2);
    \draw [connector] (sum2) -- node[name=u] {$\red{\u}$} (system);
    \node [output, right of=system, node distance=2cm] (output) {};
    \draw [connector] (input) -- node {$\green{\mathbf r}$} (sum);
    \draw [connector] (sum) -- node {$\mathbf e$}(controller);
    \draw [connector] (system) -- node [name=y] {$\red{\mathbf y}$}(output);
    \draw [connector] (output)+(-0.3cm,0) -- ++(-0.3cm,-1cm) -| node [near start] {} (sum.south);
    \draw [connector] (sum)+(0.3cm,0) |- (err) -| (sum2);
\end{tikzpicture}
}
\caption{Stability analysis of the closed-loop obtained with a reduced-order controller $\mathbf{K}$: the closed-loop is reformulated according to the controller modelling error $\mathbf{K}-\mathbf{K}^\star$.}
\label{fig:small_gain}
\end{figure}

One of the major challenges in this data-driven control strategy is to preserve internal stability. While the ideal controller is known to be stabilising thanks to the choice of an achievable reference model, see \cite{kergus2019}, there is no guarantee regarding the reduced-order controllers. Therefore, it is necessary to analyse the internal stability during the controller reduction step. To that extent, the resulting closed-loop is written as on Figure \ref{fig:small_gain}. This scheme makes the controller error appear as a perturbation. 

\begin{figure}[h]
\centering
\includegraphics[width=0.9\textwidth]{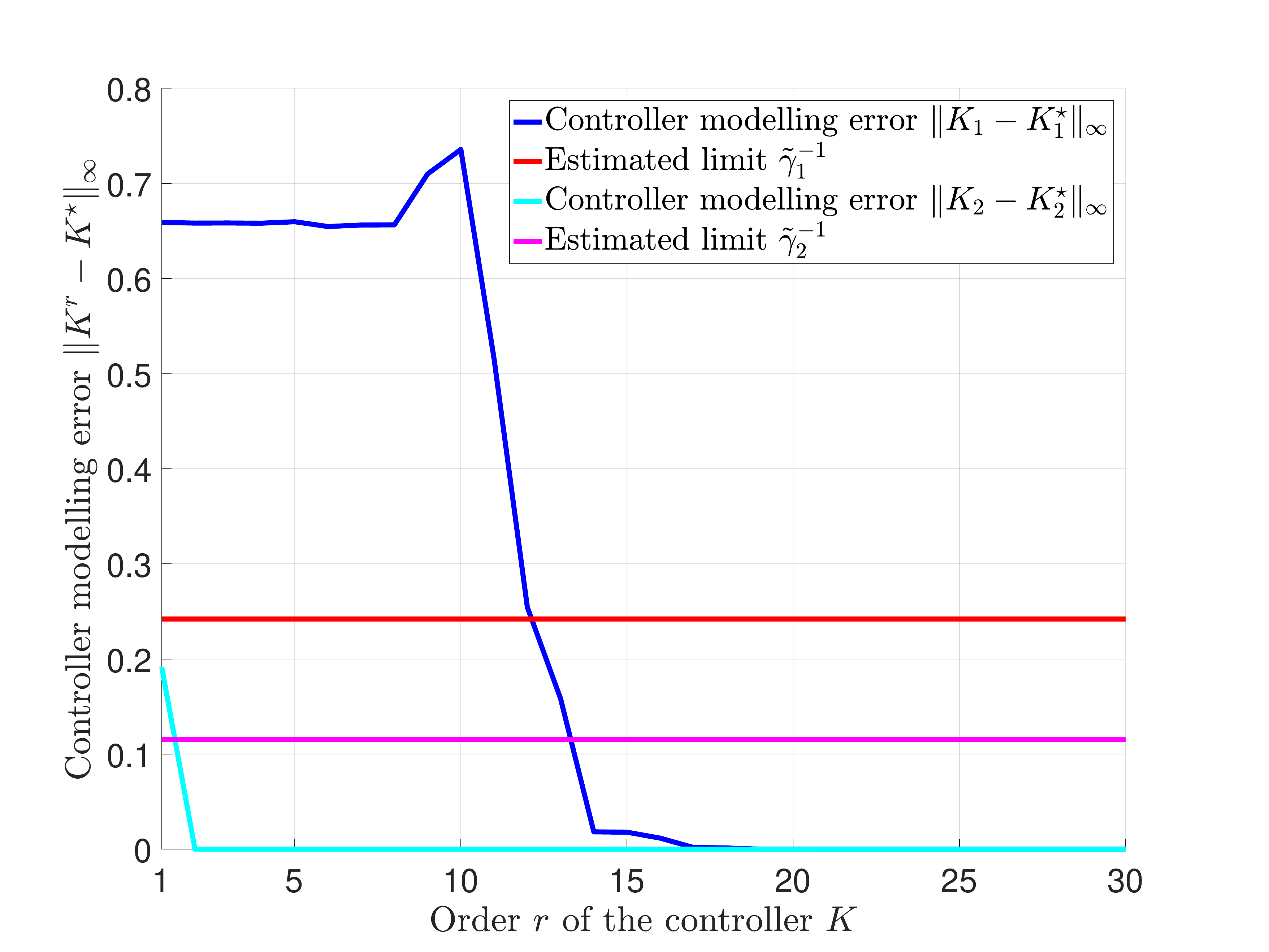}
\caption{Evolution of the controller modelling error as a function of the controller reduction step.}
\label{fig:reduction_K}
\end{figure}

It is then possible to apply the small-gain theorem: the interconnected system shown on Figure \ref{fig:small_gain} is well-posed and internally stable for all stable $\mathbf{\Delta}=\mathbf{K}-\mathbf{K}^\star$ with $\left\Vert \mathbf{\Delta} \right\Vert_\infty < \frac{1}{\gamma}$ if and only if $\left\Vert \Htran_{x_m}(1-\mathbf{M}) \right\Vert_\infty \leq\gamma$. The bound $\gamma$ on the controller modelling error can then be estimated using the data only as:
\begin{eq}
    \Tilde{\gamma}=\underset{i=1\dots m}{\textnormal{max}}\left| \Htran_{x_m}(\imath\omega_i)(1-\mathbf{M}(\imath\omega_i))\right|.
    \label{estimation_gamma}
\end{eq}
The evolution of the controller modelling error $\left\Vert \mathbf{K}-\mathbf{K}^\star \right\Vert_\infty$ according to the order of the reduced controller $\mathbf{K}$ is represented on Figure \ref{fig:reduction_K} for the two considered reference models, $\mathbf{M}_1$ and $\mathbf{M}_2$. According to Figure \ref{fig:reduction_K}, to ensure that the resulting closed-loop is internally stable, one may reduce $\mathbf{K}_1^\star$ up to $r=13$ and $\mathbf{K}_2^\star$ up to $r=2$.

However, this stability test is too conservative: it is not possible to conclude anything regarding internal stability for controllers $\mathbf{K}$ for which $\left\Vert \mathbf{K}-\mathbf{K}^\star \right\Vert_\infty>\Tilde{\gamma}$. In the present case, both order one controllers given in \eqref{LDDC_controllers} stabilise internally the rational model. Another solution would be to use the model-free stability analysis presented in \cite{cooman2018model} but it may be complicated for the user to conclude regarding internal stability. To that extent, another approach to analyse stability in a data-driven framework is introduced in Section \ref{sec:stab}. 

\subsection{Model- vs. data-driven control design remarks}

Throughout this section, it has been demonstrated how central the Loewner tool can be, either for model-driven and data-driven control. Interestingly, by choosing the closed-loop performances $\mathbf{\hat M}=\mathbf M_2$ obtained with the first approach, the second controller $\mathbf K_2$ is able to recover exactly the same controller properties, while avoiding the complex model construction step. This property reduces the time consuming model construction step and allows a quick design of the controller. However, this main advantage is balanced by the fact that in the model-based approach, the stability assessment is usually carried out using the approximated model, here $\Htranr_{x_m}^{33}$. This latter being very accurate, the eigenvalues computation is traditionally enough for concluding of the stability, robustness... On the contrary, in the second data-driven approach, the stability cannot be analysed as is and robustness (conservative) bounds are used instead.

Both approaches can be viewed as equivalent since they lead to the same  controller. Moreover, in both cases, the interpolatory framework offered by the Loewner matrices is the major ingredient for the success of the design. One may consider these approaches as complementary: the model-based approach may be privileged for critical systems where model understanding is of major  importance and for which engineering time can be spent, while the data-driven one should be the best solution for fast computation, preliminary design, for which neither safety nor critical issues are in the scope.

\section{Stability assessment of $\Ltwo$ meromorphic functions}
\label{sec:stab}

Independently of the chosen control design approach, both methodology rely on a rational model and the stability and performances obtained by the controller on the irrational model cannot be guaranteed. This is why, in this section, the stability involving the original irrational transfer $\Htran_{x_m}$, is addressed in an nonstandard manner and involving once again the Loewner matrices.

\subsection{The  \mfsa procedure}

\begin{figure}[h]
\centering
\scalebox{.8}{
\tikzstyle{block} = [draw, thick,fill=bleuONERA!20, rectangle, minimum height=4em, minimum width=6em,rounded corners]
\tikzstyle{sum} = [draw, thick,fill=bleuONERA!20, circle, node distance=1cm]
\tikzstyle{input} = [coordinate]
\tikzstyle{output} = [coordinate]
\tikzstyle{pinstyle} = [pin edge={to-,thick,black}]
\tikzstyle{connector} = [->,thick]
\begin{tikzpicture}[auto, node distance=2cm,>=latex']
    \node [input, name=input] {};
    \node [sum, right of=input] (sum) {};
    \node [block, right of=sum] (controller) {$\mathbf K(s)$};
    \node [block, right of=controller, node distance=3cm] (system) {$\Htran_{x_m}(s)$};
    \node [block, below of=controller, node distance=2cm] (delay) {\red{$e^{-\tau s}$}};
    \draw [connector] (controller) -- node[name=u] {$\red{\u}$} (system);
    \node [output, right of=system] (output) {};
    \draw [connector] (input) -- node {$\green{\mathbf r}$} (sum);
    \draw [connector] (sum) -- node {$\mathbf e$} (controller);
    \draw [connector] (system) -- node [name=y] {$\red{\y}$}(output);
    \draw [connector] (output)+(-0.3cm,0) |- node [near start] {} (delay.east);
    \draw [connector] (delay.west) -| node [near start] {} (sum.south);
\end{tikzpicture}
}
\caption{Feedback loop scheme for stability and margin computation. The objective is to assess the stability of the interconnection $\mathbf K(s)$ using $\Htran_{x_m}(s)$ and a fixed delay $\tau$ in the loop.}
\label{fig:feedbackDelay}
\end{figure}
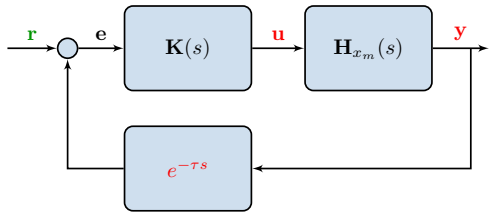

Being given a rational controller $\mathbf K$ and the irrational model defined by the meromorphic function $\Htran_{x_m}$, one important challenge is to assess the stability of the closed-loop and \eg to evaluate its stability when delay enters in the loop (delay margin of the interconnection). These questions are gathered on Figure \ref{fig:feedbackDelay}, for which the corresponding closed-loop naturally reads
\begin{equation}
\mathbf M_\tau(s) = \frac{\Htran_{x_m}(s)\mathbf K(s)}{1+\Htran_{x_m}(s)\mathbf K(s)e^{-\tau s}},
\label{eq:closedLoop}
\end{equation}
where $\tau \in \Real_+$ is the delay value affecting the loop. On the basis of \cite{PontesECC:2015} and on Chapter 5 of \cite{PoussotHDR:2019}, let us now propose a numerical procedure for the stability approximation of infinite $\Linf$ meromorphic functions. Note that the proposed version extends the one presented first in \cite{PontesECC:2015} and second in \cite{PoussotHDR:2019} by providing a much more numerically robust version, and now considers $\Linf$ functions. The \mfsa procedure given in Algorithm \ref{algo:MFSA} is first proposed. 

\begin{algorithm}
\caption{\mfsa - $\Linf$ Meromorphic Function Stability Approximation} \label{algo:MFSA} 
\begin{algorithmic}[1]
\REQUIRE $\Htran\in\Linf$, $\{\omega_i\}_{i=1}^N \in \Real_+$, $N\in \mathbb N$ and $\epsilon\in\Real_+$ (typically twice machine precision)
\STATE Sample $\Htran$ and obtain $\{\imath\omega_i,\Htran(\imath\omega_i)\}_{i=1}^N$
\STATE Construct an exact Loewner interpolant and obtain $\Htranr^r$, ensuring interpolatory conditions 
\STATE Compute $\Htranr_+^r$, the best stable approximation of $\Htranr^r$ (\eg using \cite{Kohler:2014})
\STATE Compute the stability index as $\texttt{stabTag}=\norm{\Htranr_+-\Htranr}_{\Linf}$
\STATE If $\texttt{stabTag}<\epsilon$, then $\Htran$ is stable, otherwise, $\Htran$ is unstable
\end{algorithmic}
\end{algorithm}

Algorithm \ref{algo:MFSA} embeds a relative simple procedure, which will be shown to be actually quite effective, fast and reliable. The idea consists in exactly matching the original input-output model by a rational model $\Htranr^r$, by guaranteeing interpolatory conditions. Then, to seek for the best stable approximation $\Htranr_+^r$ of the obtained model $\Htranr^r$. The $\Linf$ distance between the interpolated $\Htranr^r$ and stable $\Htranr^r_+$ models is then computed. If this latter is smaller than a given threshold $\epsilon>0$, then we conclude that $\Htran$ is stable, and unstable otherwise. By applying the procedure to \eqref{eq:closedLoop}, including the irrational model \eqref{eq:modelFreq_xm}, for varying frozen values of delay $\tau_j$ (20 linearly space between 4.6 and 5.5 has been chosen, surrounding the delay instability margin), leads to the following results of the \texttt{stabTag}:
$\big[0$, $0.8412\times 10^{-11}$, $0.2496\times 10^{-11}$, $0.4084\times 10^{-11}$, $0$, $0$, $0$, $0$, $0$, $0$, $1.7719\times 10^{7}$, $1.6346\times 10^{6}$, $50.8770$,  $34.8265$, $26.4887$, $21.3817$,  $17.9326$, $15.4471$, $13.5709$,  $12.1046\big]$. These values indicate that the closed-loop is stable up to the destabilising delay value $\tau_j\approx 5.0737$s.

\begin{figure}
\centering
\includegraphics[width=0.9\columnwidth]{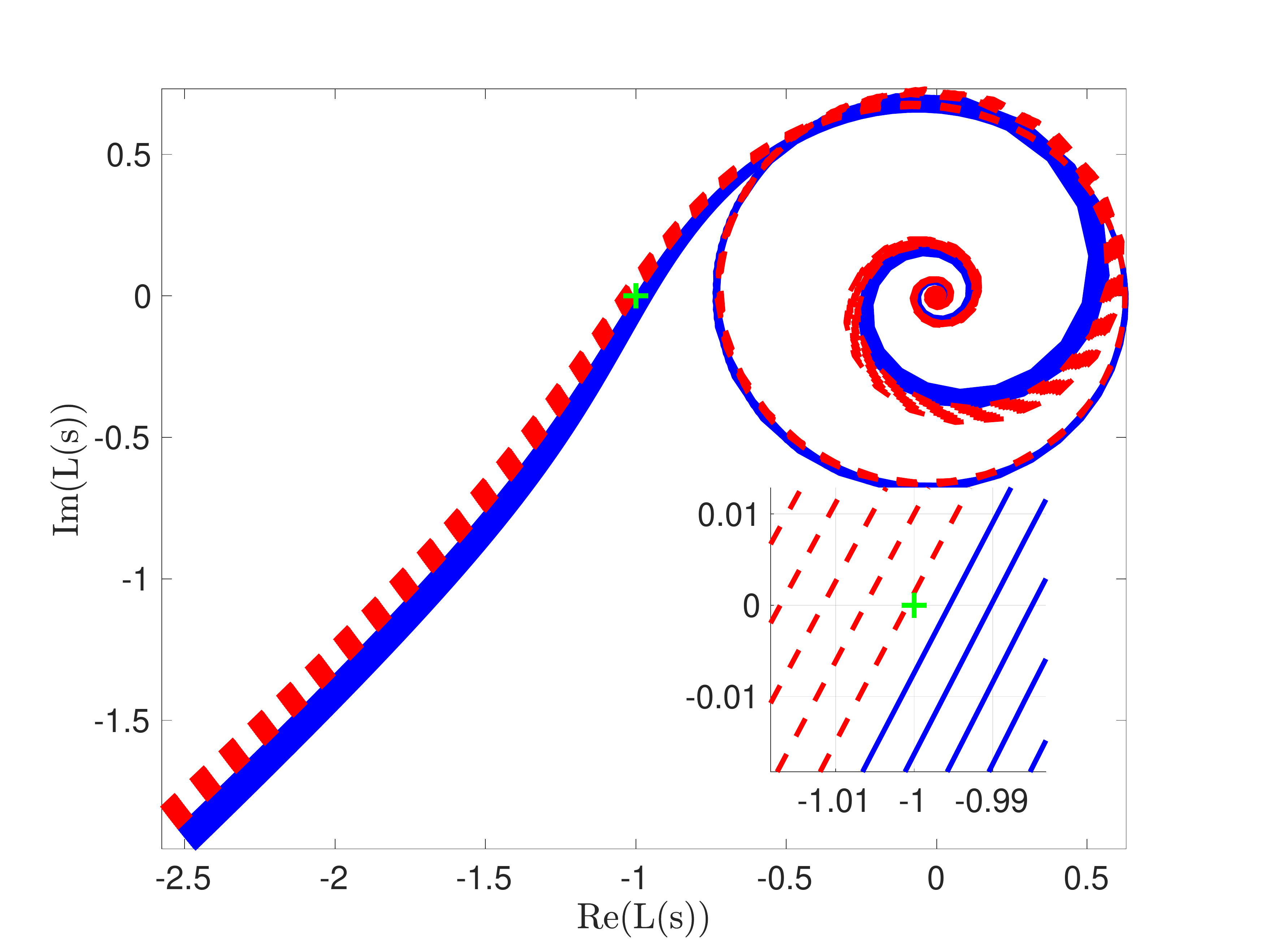}
\caption{Nyquist diagram of $\mathbf L(s)$ for varying values of $\tau$. Solid blue curves are Nyquist for which the stability tag  \texttt{stabTag} is below $10^{-10}$ (stable configuration) and dashed red curves for which tag is above $10^{-10}$. Bottom right: zoom around the Nyquist point (green $+$).}
\label{fig:nyquist}
\end{figure}

To assess this approach on such a simple \siso case, the stability can also be checked using the  Nyquist graph of $\mathbf L=\Htran_{x_m}\mathbf K$. Figure \ref{fig:nyquist}  illustrates the Nyquist curve bundle as a function of the $\tau$ and shows that the proposed approach leads to a good delay stability approximation.

\subsection{Approximation-driven arguments for the \mfsa}

\begin{figure}
\centering
\includegraphics[width=\columnwidth]{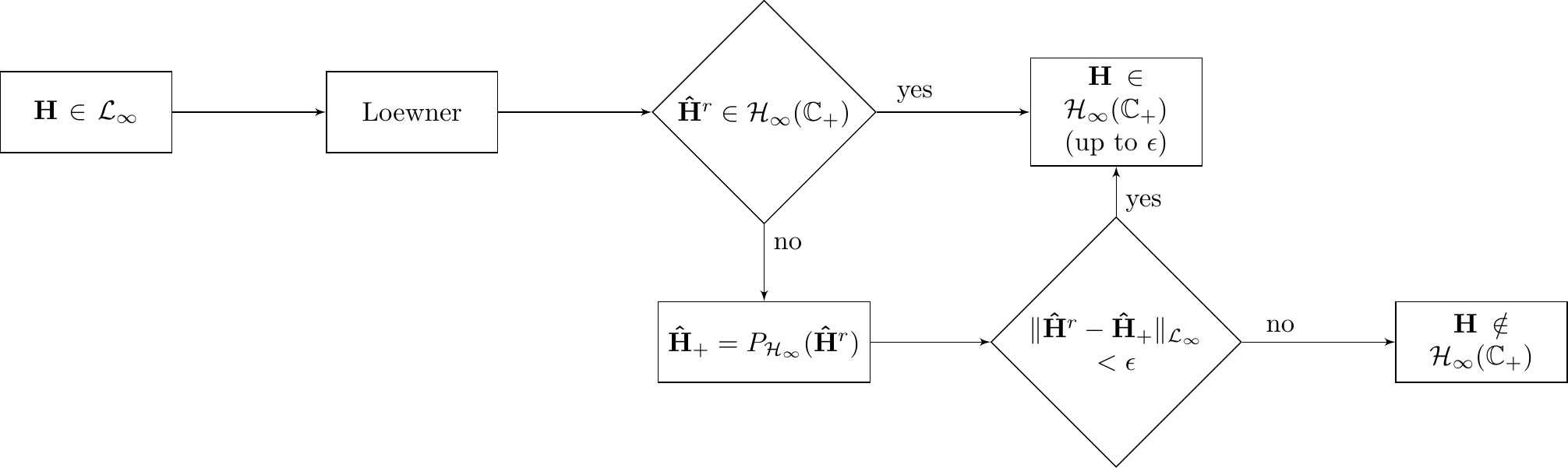}
\caption{Graphical illustration of the \mfsa procedure (starting point, at top left). }
\label{fig:process}
\end{figure}

Algorithm \ref{algo:MFSA} is rather simple and deserves some comments arguing of its viability. Let us first refer to \cite{PontesECC:2015}  where arguments and a similar procedure, involving the \tfirka algorithm \cite{BeattieCDC:2012} have been suggested (this procedure ensures $\Htwo$-optimal bi-tangential interpolatory conditions). This later provides good results but lacks in determining the approximation order $r$. Additionally, as \tfirka is an $\Htwo$-oriented procedure its validity in the $\Ltwo(\imath\Real)$ function space is limited. 
Here, the path presented in Figure \ref{fig:process} is used to construct the \mfsa procedure.

With reference to Figure \ref{fig:process} (top left three blocks), the \tfirka interpolatory conditions are released and the Loewner framework is used instead. First, one major benefit of such a trade stands in the automatic selection of the approximating order $r$, done by a rank revealing factorisation where machine precision is expected. Second, it  interpolates the data without any stability constraint. One obtains $\Htranr^r$ which tangentially interpolates the data and which, by increasing the numbers of samples $z_i$ in \eqref{eq:data} hopefully converges to the same model $\Htran^r$ (\eg the rational approximation obtained by Loewner matrices is not affected by the number of interpolation points). At this point, the resulting model $\Htranr^r$ may be stable or unstable. 
\begin{itemize}
\item If $\Htranr^r$ is stable, then one concludes on the stability of $\Htran$, up to some tolerance $\epsilon >0$. Indeed, one could always add $\epsilon \mathbf{H}_a \in \Hinf$ to the data without being able to see the effects on on the interpolation conditions achieved by $\Htranr^r$, due to numerical accuracy.
\item If $\Htranr^r$ is unstable, as it is always possible to find an unstable approximant to a stable model in the $\Ltwo$ sense \cite{PontesECC:2015}, one may use the projection onto a stable subspace, here the $\mathcal{H}_\infty$ one, denoted $\Htranr_+ = P_\Hinf(\Htranr^r)$, to emphasise the importance of the unstable part of $\Htranr^r$ on the interpolation conditions. If the unstable part plays a negligible role (\ie $\norm{\Htranr^r-\Htranr_+}_\Linf < \epsilon$), then a stable interpolating model has been found for $\Htran$ which is therefore likely to be stable. Otherwise, if the unstable part cannot be removed, then it is likely that $\Htran$ is unstable.

\end{itemize}

The Loewner framework allows to find a rational model $\Htranr\in\mathcal{RL}_\infty$ that interpolates $\Htran\in\Linf$ at an arbitrary number of frequencies. The suggestion we claim is in twofold. One is always able to find a rational model $\Htranr^r \in \mathcal{RL}_\infty$ that well reproduces $\Htran \in \Linf$(at least interpolates a large number of points). We assume that this implies a convergence in the $\Linf$ sense, meaning that one is always able to find a rational function matching an irrational one defined over $\Linf$. Then, if $\Htranr^r \in \mathcal{RL}_\infty$ can in addition be projected onto $\Htranr_+=P_\Hinf(\Htranr^r) \in\mathcal{H}_\infty(\Cplx_+)$ with negligible loss in the $\Linf$-norm, then the unstable part can be considered as irrelevant for the behaviour description of $\Htran$, which can thus be assumed stable. Otherwise, if the unstable part leads to important $\Linf$-norm mismatch, $\Htran$ can be considered as unstable. 




\section{Conclusion}

In this chapter, the interpolatory framework proposed by the Loewner setup, as introduced in the seminal paper \cite{Mayo:2007}, has been further used for the control design and the stability estimation. The singularity of the proposed approach is to show that the Loewner is not only a model approximation tool, but a complete dynamical-oriented tool. The main contributions of this work is twofold. First, to compare frequency-oriented data- and model-based control design approaches, showing that both lead to similar performances. Second, to suggest a method for approximating the stability of any $\Linf$ functions (either rational or irrational). Both contributions are based on Loewner matrices. Through an academic example described by a linear \pde set, the Loewner framework has been used for different purpose, showing its impressive versatility and applicability to solve complex problems. To the authors perspective, this approach opens the fields for analysing and controlling irrational (infinite-dimensional) models in a relatively simple manner. Even if the approach does not stand as a completely closed solution, it can be viewed as an alternative for engineers and practitioners to deal with irrational models in a simple manner.

\bibliographystyle{plain}

\end{document}